\documentclass[12pt]{amsart}
\usepackage{amssymb}
\usepackage{latexsym}
\usepackage{stix}

\usepackage{color}

\oddsidemargin=\evensidemargin
\begin{document}
\newfont{\blb}{msbm10 scaled\magstep1} 

\newtheorem{defi}{Definition}
\newtheorem{theo}{Theorem}[section]
\newtheorem{theo1}{Theorem}[section]
\newtheorem{prop}[theo]{Proposition}
\newtheorem{lemm}[theo]{Lemma}
\newtheorem{coro}[theo]{Corollary}
\date{}
\author[I. de las Heras and G. Traustason]{Iker de las Heras \\ University of the Basque Country, Spain \\  \\ Gunnar Traustason \\ University of Bath, UK}

\title{Powerfully solvable and powerfully simple groups}

\newcommand{\comment}[1]{{\color{red}\rule[-0.5ex]{3pt}{2.5ex}}\marginpar{\small\begin{flushleft}\color{red}#1\end{flushleft}}}

\setlength{\marginparwidth}{2cm}

\newcommand{\F}{\mathbb{F}}

\thanks{The first author is supported by the Spanish Government grant MTM2017-86802-P, partly with FEDER funds, and by the Basque Government grant IT974-16.
He is also supported by a predoctoral grant of the University of the Basque Country}

\maketitle
\begin{abstract}
We introduce the notion of a powerfully solvable group. These are powerful groups possessing an abelian series of a special kind. These groups include in particular the class of powerfully nilpotent groups. We will also see that for a certain rich class of powerful groups we can naturally introduce the term powerfully simple group and prove a Jordan-H\"{o}lder type theorem that justifies the term.
\end{abstract}
\section{Introduction}
In this paper $p$ is always an odd prime. \\ \\
Recall [2] that a finite $p$-group is {\it powerful} if $[G,G]\leq G^{p}$. More generally, a subgroup $H$ of $G$ is powerfully embedded in $G$ if $[H,G]\leq H^{p}$.  \\ \\
The following useful property, known as Shalev's Interchange Lemma [3], will be used a number of times in this paper: If $H$ and $K$ are powerfully embedded in $G$ then $[H^{p^{i}},K^{p^{j}}]=[H,K]^{p^{i+j}}$. \\ \\
Another term that we need is a {\it powerful basis}. For any powerful $p$-group $G$, there exist generators $a_{1},\ldots ,a_{r}$ 
such that $G=\langle a_{1}\rangle\cdots \langle a_{r}\rangle$ and where $|G|=o(a_{1})\cdots o(a_{r})$. The number of generators of any given order is an invariant. We say that $G$ is of {\it type} $(1,\stackrel{r_{1}}{\ldots},1,2,\stackrel{r_{2}}{\ldots},2,\ldots )$
if there are $r_{i}$ generators of order $p^{i}$. \\ \\
In  [5] the notion of powerful nilpotence and powerfully central chain was introduced. If $K\leq H\leq G$, then a chain of subgroups
                          $$H=H_{0}\geq H_{1}\geq \cdots \geq H_{n}=K$$
is {\it powerfully central in $G$}  if $[H_{i},G]\leq H_{i+1}^{p}$ for $i=0,\ldots ,n-1$. A finite $p$-group is said to be {\it powerfully nilpotent} if it has a powerfully central series $G=G_{0}\geq G_{1}\geq \cdots \geq G_{n}=1$. The smallest possible length of such a series is then called the {\it powerful nilpotence class} of $G$. \\ \\
Powerful nilpotence  leads then naturally to a classification in terms of an ancestry tree and powerful co-class. In [5] it was shown that for every prime $p$ there are finitely many powerfully nilpotent $p$-groups of each powerful co-class, and some general theory was developed for powerfully nilpotent groups.
In [6] the powerfully nilpotent groups of maximal powerful class are studied. \\ \\
In this paper we consider a natural larger class of powerful $p$-groups. \\ \\
{\bf Definition}. Let $G$ be a finite $p$-group and $K\leq H\leq G$. We say that a chain 
                  $$H=H_{0}\geq H_{1}\geq \cdots \geq H_{n}=K$$
is {\it powerfully abelian} if $[H_{i},H_{i}]=H_{i+1}^{p}$ for $i=0,\ldots ,n-1$. \\ \\
{\bf Definition}. A finite $p$-group $G$ is {\it powerfully solvable}, if there exists a powerfully abelian chain
      $$G=G_{0}\geq G_{1}\geq\cdots \geq G_{n}=1.$$
The smallest possible length $n$ is called the {\it powerful derived length} of $G$. \\ \\
The structure of the paper is as follows. In Section 2 we show that all powerful $p$-groups of rank $2$ are powerfully solvable and based on the work in [7] we provide a classification of all these groups as well as a closed formula for the number of
such groups of order $p^{x}$. In Section 3 we introduce the notion of a powerfully solvable presentation that will be useful later on when going through some classification and
calculating growth. In Section 4 we classify all powerful groups  of order up to $p^{5}$ and see that these are all powerfully solvable. In Section 5 we discuss the growth of powerfully solvable  groups, and various other classes of powerful groups, that are
of exponent $p^{2}$. In Section 6 we consider the rich class ${\mathcal P}$ of all powerful $p$ groups of type $(2,2,\ldots ,2)$  and see that powerful nilpotence and powerful solvability play a similar role here as nilpotence and solvability for the class of all groups. The notion of a powerfully simple group arises naturally and we are able to prove a Jordan-H\"{o}lder like result that justifies the term. Finally in Section 7 we classify all the powerfully simple groups of order $p^{6}$. The number turns out to depend on the prime $p$. 
\section{Powerful groups of rank $2$}
It turns out that all powerful $p$-groups of rank $2$ are powerfully solvable. In fact something stronger is true.
\begin{prop} Let $G$ be a powerful $p$-group. If $[G,G]$ is cyclic then $G$ is powerfully solvable of powerful derived length at most
$2$. 
\end{prop} 
\noindent{\bf Proof}\ \ As $G$ is powerful we have that $[G,G]=\langle g^{p}\rangle$ for some $g\in G$. Therefore
         $$G\geq \langle g\rangle\geq 1$$
is a powerfully abelian chain. $\Box$ \\ \\
In [7] the powerfully nilpotent groups of rank $2$ are classified and a closed formula is given for the number 
of powerfully nilpotent groups of order $p^{x}$. In fact there is implicitly the following  classification of all powerful $p$-groups of
rank $2$. By Proposition 2.1 we know that these are all powerfully solvable. \\ \\
{\bf Classification of the non-abelian powerful groups of rank $2$}.  \\ \\
(I) Semidirect products: 
                       $$G=\langle a,b:\,a^{p^{n}}=b^{p^{m}}=1,\,[a,b]=a^{p^{r}}\rangle $$
with $n-r\leq m$ and $1\leq r\leq n-1$. \\ \\
(II) Non-semidirect products: 
               $$G=\langle a,b:\,a^{p^{n}}=1,\,b^{p^{m}}=a^{p^{l}},\, [a,b]=a^{p^{r}}\rangle$$
with $1\leq r<l\leq n-1$ and $n-r\leq l < m$. \\ \\ \\
From [7] we also know that a group above is powerfully nilpotent if and only if $r\geq 2$. Thus, it is easy to determine that 
there are $\lfloor\frac{x-1}{2}\rfloor$ semidirect products and $\lfloor\frac{x-4}{2}\rfloor$ non-semidirect products of order $p^{x}$ that are not powerfully nilpotent (here $\lfloor\cdot\rfloor$ stands for the floor function). From this, the discussion above and  [7, Proposition 2.2] we also get the following. \\ \\
{\bf Enumeration}.  For $x\geq 3$, the number of powerful $p$-groups of rank $2$ and order $p^{x}$ is 
\begin{eqnarray*}
       \frac{x^{3}+12x^{2}+12x}{72} & \mbox{if} & x\equiv 0\pmod{6}\\
       \frac{x^{3}+12x^{2}+3x-16}{72} & \mbox{if} & x\equiv 1\pmod{6} \\
      \frac{x^{3}+12x^{2}+12x-8}{72} & \mbox{if} & x\equiv 2\pmod{6} \\
     \frac{x^{3}+12x^{2}+3x}{72} & \mbox{if} & x\equiv 3 \pmod{6} \\
     \frac{x^{3}+12x^{2}+12x-16}{72} & \mbox{if} & x\equiv 4 \pmod{6} \\
      \frac{x^{3}+12x^{2}+3x-8}{72} & \mbox{if} & x\equiv 5\pmod{6}. \\
\end{eqnarray*}

\section{Presentations}
\begin{lemm} Let $G$ be a finite $p$-group and let $K<H\leq G$ where $[H,H]\leq K^{p}$. If for some positive integer $n$ we have $K^{p^{n}}=H^{p^{n}}$, then there exists $x\in K\setminus H$ such that $x^{p^{n}}=1$.
\end{lemm}
\noindent{\bf Proof}\ \ We prove this by induction on $n$. Suppose first that $n=1$. We will  show that for every $j\geq 1$, there exists $x\in K\setminus H$ such that $x^{p}\in H^{p^{j}}$. For $j=1$ this is immediate from the hypothesis, so assume by induction on $j\geq 2$ that 
we know that $x^{p}\in H^{p^{j-1}}$ for some $x\in K\setminus H$. Then there exists $y\in H^{p^{j-2}}$ such that $x^{p}=y^{p}$. Now, by the Hall-Petresco Identity, as $p>2$, we have 
                      $$(xy^{-1})^{p}=x^{p}y^{-p}c_{2}^{p\choose 2}c_{3}^{p\choose 3}\cdots c_{p}^{p\choose p},$$
where $c_{k}\in [H^{p^{j-2}},K,\stackrel{k-1}{\ldots},K]$ for $k=2,\ldots ,p$. Notice that $K$ is powerful and that $H$ is powerfully embedded in $K$. We can thus use  Shalev's Interchange Lemma. Therefore for $2\leq k\leq p-1$ we have
                            $$c_{k}^{p\choose k}\in [H^{p^{j-2}},K]^{p}=[H,K]^{p^{j-1}}\leq H^{p^{j}}$$
and 
                          $$c_{p}\in [H^{p^{j-2}},K,K]=[H,K,K]^{p^{j-2}}\leq [H^{p},K]^{p^{j-2}}=[H,K]^{p^{j-1}}\leq H^{p^{j}}.$$
Since $x^{p}y^{-p}=1$, we then have $(xy^{-1})^{p}\in H^{p^{j}}$. This finishes the inductive step. Taking $j$ such that $H^{p^{j}}=1$ we see that we can pick $x\in K\setminus H$ such that $x^{p}=1$. \\ \\
Now suppose $n>1$ and that the result holds for smaller values of $n$. If $K^{p^{n-1}}=H^{p^{n-1}}$, then by induction hypothesis we know there exists $x\in K\setminus H$ where $x^{p^{n-1}}=1$ and thus $x^{p^{n}}=1$. We can thus assume that $K^{p^{n-1}}\not =H^{p^{n-1}}$. Now
$$[K^{p^{n-1}},K^{p^{n-1}}]\leq [K,K]^{p^{2n-2}}\leq H^{p^{2n-1}}\leq (H^{p^{n-1}})^{p}.$$
Therefore by the induction hypothesis there exists an element $y\in K^{p^{n-1}}\setminus H^{p^{n-1}}$
such that $y^{p}=1$. Since $K$ is powerful we have $y=x^{p^{n-1}}$ for some $x\in K\setminus H$ and then $x^{p^{n}}=y^{p}=1$. $\Box$
\begin{theo}
\label{theorem presentation}
Let $G$ be a finite  $p$-group of rank $r$ and exponent $p^{e}$ where $G/G^{p^{2}}$ is powerfully solvable. Then $G$ is powerfully solvable. Furthermore, we can choose our generators $a_{1},a_{2},\ldots ,a_{r}$ such that $|G|=o(a_{1})\cdots o(a_{r})$ and such that the chain
 $$\begin{array}{lllll}
       G=\langle a_{1},a_{2},\ldots ,a_{r}\rangle
       &
       \geq
       &
       \langle a_{1}^{p},a_{2},\ldots ,a_{r}\rangle
       &
       \geq \cdots \geq  & G^{p} \\
      G^{p}=\langle a_{1}^{p},a_{2}^{p},\ldots ,a_{r}^{p}\rangle
      &
      \geq
      &
      \langle a_{1}^{p^{2}},a_{2}^{p},\ldots ,a_{r}^{p}\rangle
      &
      \geq \cdots \geq  & G^{p^{2}} \\
     &   &   \vdots &  & \\
    G^{p^{e-1}}=\langle a_{1}^{p^{e-1}},a_{2}^{p^{e-1}},\ldots ,a_{r}^{p^{e-1}}\rangle
    &
    \geq
    &
    \langle a_{1}^{p^{e}},a_{2}^{p^{e-1}},\ldots ,a_{r}^{p^{e-1}}\rangle
    &
    \geq \cdots \geq  & G^{p^{e}}=1 
\end{array}$$
is powerfully abelian.
\end{theo}
\noindent{\bf Proof}\ \ Suppose, using the fact that $G/G^{p^{2}}$ is powerfully solvable, that $G=K_{0}>K_{1}>\cdots >K_{m}=G^{p^{2}}$ is a chain that is powerfully abelian modulo $G^{p^{2}}$. Notice that $[G,G]\leq K_{1}^{p}G^{p^{2}}\leq G^{p}$ and the group is thus  powerful. 
 In particular, we have $[G^{p},G]\leq G^{p^{2}}$ and $(G^{p})^{p}=G^{p^{2}}$. Therefore $G=K_{0}G^{p}\geq K_{1}G^{p}\geq \cdots \geq K_{m}G^{p}=G^{p}$ is also powerfully abelian. Removing redundant terms and refining if necessary, we get a powerfully abelian chain
                      $$G=H_{0}>H_{1}>\cdots >H_{r}=G^{p}$$
where the factors are of size $p$. Now notice that for $0\leq i\leq r-1$ and $0\leq j\leq e$ we have $[H_{i}^{p^{j}},H_{i}^{p^{j}}]=
[H_{i},H_{i}]^{p^{2j}}\leq H_{i+1}^{p^{j+1}}$. This gives us the powerfully abelian chain we wanted. It remains to see that we can furthermore pick our generators such that $a_{1},\ldots ,a_{r}$ is a powerful basis for $G$. Let us pick our generators of $G$ such that  for every $1\le i\le r-1$ we have $H_{i}=\langle a_{i+1},\ldots ,a_{r}\rangle G^{p}$.   If $H_{i}^{p}=H_{i+1}^{p}$ for some $1\leq i\leq r-1$ then we can know from Lemma 3.1 that we can pick $a_{i+1}$ such that $a_{i+1}^{p}=1$. We can also
in that case move the generator in front. We thus have 
              $$a_{1}^{p}=\cdots =a_{r-r_{1}}^{p}=1\mbox{ and }H_{r-r_{1}}^{p}>\cdots >H_{r}^{p}=G^{p},$$
where $r_{1}=\mbox{rank\,}(G^{p})$. Now consider the chain
          $$G^{p^{2}}=H_{r-r_{1}}^{p^{2}}\geq \cdots \geq H_{r}^{p^{2}}=G^{p^{3}}.$$
Again if $H_{i}^{p^{2}}=H_{i+1}^{p^{2}}$, then we know by Lemma 3.1 that we can pick $a_{i+1}$ such that $a_{i+1}^{p^{2}}=1$.
Continuing in this manner we see that we can choose our generators such that for $1\leq i\leq r$ we have $o(a_{i})=p^{j}$ where $j$ is the smallest positive integer such that $H_{i-1}^{p^{j}}=H_{i}^{p^{j}}$. Also we have that $\mbox{rank\,}(G^{p^{i}})$ is then the number of $1\leq i\leq r$ such that $a_{i}^{p^{j}}\not =1$. Let $r_{j}$ be the rank of $G^{p^{j}}$. Then ($r_{e}=0$)
    $$|G|=p^{r_{0}+r_{1}+\cdots +r_{e-1}}=p^{r_{0}-r_{1}}\cdot (p^{2})^{r_{1}-r_{2}}\cdots (p^{e-1})^{r_{e-1}-r_{e}}
             =o(a_{1})\cdots o(a_{r}).$$
This finishes the proof. $\Box$ \\ \\
{\bf Powerfully solvable presentations}. It follows in particular from Theorem  3.2 that a powerfully solvable group
of order $p^{n}$ and rank $r$ has a presentation with generators $a_{1},\ldots , a_{r}$ and relations
\begin{equation}
                    a_{1}^{p^{n_{1}}}=1,\ldots ,a_{r}^{p^{n_{r}}}=1
\end{equation}
and 
\begin{equation}
                   [a_{j},a_{i}]=a_{1}^{m_{1}(i,j)}\cdots a_{r}^{m_{r}(i,j)},\ 1\leq i<j\leq r,
\end{equation}
where all the power indices $m_{k}(i,j)$ are divisible by $p$ and where furthermore $p^{2}|m_{k}(i,j)$ whenever $k\leq i$. Notice that  $G$ is the largest finite $p$-group satisfying these relations. To see this let $H$ be the largest finite $p$-group satisfying these relations. The group $H/H^{p^{2}}$ is powerfully solvable and thus $H$ is powerfully solvable by  Theorem 3.2. In particular $H$ is powerful and therefore $|H|\leq o(a_{1})\cdots o(a_{r})$. However $G$ is a homomorphic image of $H$ and thus $|H|=o(a_{1})\cdots o(a_{r})$. Hence $H$ is isomorphic to $G$.  \\ \\
A presentation with generators $a_{1},\ldots ,a_{r}$ and relations of the form (1) and (2) is called a {\it powerfully solvable presentation}. We say that such a presentation is \emph{consistent} if the presentation determines a group of order $p^{n_{1}}\cdots p^{n_{r}}$. 
\section{Classification of powerful  groups of order up to $p^{5}$}
In this section we will find all powerful $p$-groups of order up to and including $p^{5}$. It turns out that these are all powerfully
solvable. We will see later that there are many powerful groups of order $p^{6}$ that are not powerfully solvable. Let us now
turn to our task in this section. There are 2 non-abelian groups of order $p^{3}$. The Heisenberg group of exponent $p$ cannot be powerful as it is of exponent $p$. The other group is a semidirect product of a cyclic group of order $p^{2}$ by a cyclic group
of order $p$:
                                     $$G_{1}=\langle a,b:\,a^{p^{2}}=b^{p}=1,\ [a,b]=a^{p}\rangle.$$
Notice that this group is powerfully solvable with a powerfully abelian chain $G>\langle a\rangle>1$. It is however not powerfully
nilpotent as $Z(G)^{p}=1$. Adding to this the $3$ abelian groups of order $p^{3}$, we see that there are in total 4 powerfully solvable groups of order $p^{3}$. \\ \\
Before moving on we consider a general setting like in [7] that includes a number of groups that will occur, namely the non-abelian groups of type 
$(1,\stackrel{t}{\ldots},1,n)$ where $n$ is an integer greater than $1$. Suppose 
         $$G=\langle a_{1},\ldots ,a_{t},b\rangle$$
is a powerful group of this type where $a_{i}$ is of order $p$ and $b$ of order $p^{n}$. Notice that $G^{p}=\langle b^{p}\rangle$ is cyclic and it follows from [5, Corollary 3.3] that $G^{p}\leq Z(G)$.  In particular $G$ is nilpotent of class at most $2$ and $[G,G]^{p}=[G^{p},G]=1$. Next observe that $\Omega_{1}(G)=\langle a_{1},\ldots ,a_{r},b^{p^{n-1}}\rangle$ where 
$\Omega_{1}(G)$ is the subgroup consisting of all elements of order dividing $p$. Thus $[G,G]=\langle b^{p^{n-1}}\rangle$. Now, consider the vector space $V=\Omega_{1}(G)G^{p}/G^{p}$ over the field $\F_p$ of $p$ elements. The commutator operation naturally induces an alternating form on $V$ through 
      $$(xG^{p},yG^{p})=\lambda\mbox{\ if\ }[x,y]=b^{\lambda p^{n-1}}.$$
Without loss of generality we can suppose that our generators have been chosen such that we get the following orthogonal decomposition 
                 $$V=\langle a_{1}G^{p},a_{2}G^{p}\rangle\oplus \cdots \oplus \langle a_{2s-1}G^{p},a_{2s}G^{p}\rangle\oplus V^{\perp}$$
where $V^{\perp}=\langle a_{2s+1}G^{p},\ldots ,a_{t}G^{p}\rangle$ and $(a_{2i-1}G^{p},a_{2i}G^{p})=1$ for $i=1,\ldots ,s$. There are now two cases to consider, depending on whether or not $Z(G)\leq \Omega_{n-1}(G)$. \\ \\
Suppose first that $Z(G)\not\leq\Omega_{n-1}(G)$. This means that $Z(G)$ contains some element $b^{l}u$ with $u\in \langle a_{1},\ldots ,a_{r}\rangle$ and $0<l<p$. Thus without loss of generality we can assume that $b\in Z(G)$. We thus get a powerful group $G=A(n,t,s)$ with relations 
$$\begin{array}{l}
    a_{1}^{p}=\cdots =a_{t}^{p}=b^{p^{n}}=1; \\
 \mbox{}  [a_{2i-1},a_{2i}]=b^{p^{n-1}}\mbox{\ for\ }i=1,\ldots ,s; \\
  \mbox{}[a_{i},a_{j}]=1\mbox{\ otherwise for\ }1\leq i<j\leq t; \\
\mbox{}  [a_{i},b]=1\mbox{\ for }1\leq i\leq t.
\end{array}$$
Notice that we have $\langle b\rangle\leq Z(G)$ and $[G,G]\leq \langle b^{p}\rangle$ and thus these groups are all powerfully nilpotent, as was observed in [7]. Notice that for a fixed $n\geq 2$ and $t\geq 2$ we get $\lfloor t/2\rfloor$ such groups. \\ \\
We then consider the case when $Z(G)\leq \Omega_{n-1}(G)$. Notice first that replacing $b$ by a suitable $ba_{1}^{\alpha_{1}}\cdots a_{2s}^{\alpha_{2s}}$, we can assume that $b$ commutes with $a_{1},\ldots ,a_{2s}$. As $b\not\in Z(G)$ we then must have $t>2s$ and similarly, replacing $a_i$ by a suitable $a_{2s+1}^{\alpha_{2s+1}}\cdots a_t^{\alpha_t}$, we can pick our generators $a_{2s+1},\ldots, a_{t}$ such that $[a_{2s+1},b]=b^{p^{n-1}}$ and $[a_{2s+2},b]=\cdots =[a_{t},b]=1$. We thus arrive at a group $G=B(n,t,s)$ satisfying the relations
$$\begin{array}{l}
    a_{1}^{p}=\cdots =a_{t}^{p}=b^{p^{n}}=1; \\
 \mbox{}  [a_{2i-1},a_{2i}]=b^{p^{n-1}}\mbox{\ for\ }i=1,\ldots ,s; \\
  \mbox{}[a_{i},a_{j}]=1\mbox{\ otherwise for\ }1\leq i<j\leq t; \\
\mbox{}  [a_{1},b]=\cdots =[a_{2s},b]=[a_{2s+2},b]=\cdots =[a_{t},b]=1 \mbox{\ for }1\leq i\leq t; \\
\mbox{} [a_{2s+1},b]=b^{p^{n-1}}.
\end{array}$$
Notice that for a fixed $n\geq 2$ and $t\geq 1$ there are $\lfloor (t+1)/2\rfloor$ such groups. Notice also that when $n\geq 3$ then the group is powerfully nilpotent as $\langle b^{p}\rangle\leq Z(G)$ and $[G,G]\leq \langle b^{p^{2}}\rangle$. 
For $n=2$ this is not the case but the group is still powerfully solvable as we have a powerfully abelian chain $G>\langle b\rangle>1$ with $[G,G]\leq \langle b^{p}\rangle$. We are now ready for groups of order $p^{4}$.  In the following we will omit writing relations of the form $[x,y]=1$. \\ \\
{\bf Groups of order $p^{4}$}. From our analysis of non-abelian groups of  rank $2$ we get two such groups:
     $$G_{2}=\langle a,b:\, a^{p^{2}}=b^{p^{2}}=1,\, [a,b]=a^{p}\rangle\mbox{\ \ and\  \  }G_{3}=\langle a,b:\, a^{p^{3}}=b^{p}=1,\ [a,b]=a^{p^{2}}\rangle.$$
Here $G_{3}$ is furthermore powerfully nilpotent. The only non-abelian groups apart from these are of type $(1,1,2)$ and from the analysis of such groups above we know there are two groups:
        $$G_{4}=A(2,2,1)=\langle a,b,c:\, a^{p}=b^{p}=c^{p^{2}}=1,\, [a,b]=c^{p}\rangle,$$
and 
      $$G_{5}=B(2,2,0)=\langle a,b,c:\,a^{p}=b^{p}=c^{p^{2}}=1,\,[a,c]=c^{p}\rangle .$$
Apart from these there are 5 abelian groups and we thus get in total $\boldsymbol{9}$ groups.  \\ \\
{\bf Groups of order $p^{5}$}.  Again our analysis of groups of rank $2$ and those of type $(1,1,3)$ and $(1,1,1,2)$ gives
us the following non-abelian powerfully solvable groups:
$$\begin{array}{l}G_{6}=\langle a,b:\,a^{p^{2}}=b^{p^{3}}=1,\, [a,b]=a^{p}\rangle,\ 
    G_{7}=\langle a,b:\,a^{p^{3}}=b^{p^{2}}=1,\, [a,b]=a^{p}\rangle, \\
      G_{8}=\langle a,b:\,a^{p^{3}}=b^{p^{2}}=1,\,[a,b]=a^{p^{2}}\rangle,\ G_{9}=
   \langle a,b:\,a^{p^{4}}=b^{p}=1,\ [a,b]=a^{p^{3}}\rangle,\end{array}$$
and
$$\begin{array}{l}
   G_{10}=A(3,2,1)=\langle a,b,c:\,a^{p}=b^{p}=c^{p^{3}}=1,\,[a,b]=c^{p^{2}}\rangle; \\
  G_{11}=B(3,2,0)=\langle a,b,c:\,a^{p}=b^{p}=c^{p^{3}}=1,\,[a,c]=c^{p^{2}}\rangle; \\
   G_{12}=A(2,3,1)=\langle  a,b,c,d:\,a^{p}=b^{p}=c^{p}=d^{p^{2}}=1,\,[a,b]=d^{p}\rangle; \\
  G_{13}=B(2,3,0)=\langle a,b,c,d:\,a^{p}=b^{b}=c^{p}=c^{p^{2}}=1,\, [a,b]=d^{p}, [c,d]=d^{p}\rangle; \\
G_{14}=B(2,3,1)=\langle a,b,c,d:\,a^{p}=b^{p}=c^{p}=d^{p^{2}}=1,\,[a,b]=d^{p}, [c,d]=d^{p}\rangle.
\end{array}$$
Here $G_{8},G_{9},G_{10},G_{11},G_{12}$ are furthermore powerfully nilpotent. Apart from these $9$ groups, there are $7$ abelian groups. 
We are now only left with the non-abelian groups of type $(1,2,2)$ that will contain a number of different groups and we need to deal with a number of subcases. \\ \\
Suppose that we have generators $a,b,c$ of orders $p,p^{2},p^{2}$. \\ \\
{\it Case 1}. ($Z(G)^{p}\not =1$). Notice that we then must have $|Z(G)^{p}|=p$ as otherwise $G/Z(G)$ is cyclic and thus $G$ abelian. 
We can assume that $c\in Z(G)$ and that $Z(G)^{p}=\langle c^{p}\rangle$. Notice also that $[G,G]=\langle [a,b]\rangle$ is cyclic. There are two possibilities. On the one hand, if $[G,G]\leq Z(G)^{p}$, then we can choose our generators so that we get 
a group with the following presentation:
                              $$G_{15}=\langle a,b,c:\, a^{p}=b^{p^{2}}=c^{p^{2}}=1,\, [a,b]=c^{p}\rangle.$$
On the other hand if $[G,G]\not\leq Z(G)^{p}$, it is not difficult to see that we can pick our generators so that we get a group with the presentation
                        $$G_{16}=\langle a,b,c:\, a^{p}=b^{p^{2}}=c^{p^{2}}=1,\, [a,b]=b^{p}\rangle.$$
Notice that both these groups are powerfully solvable and that $G_{15}$ is furthermore powerfully nilpotent. \\ \\
{\it Case 2}. ($Z(G)^{p}=1$ and $G/Z(G)$ has rank $2$). Then we must have $a\in Z(G)$. It is not difficult to see that
in this case we can choose $b,c$ such that $[b,c]=c^{p}$ and we get the powerfully solvable group 
                     $$G_{17}=\langle a,b,c:\,a^{p}=b^{p^{2}}=c^{p^{2}}=1,\, [b,c]=c^{p}\rangle.$$  
Before considering further cases, we first show that if $Z(G)^{p}=1$ and $G/Z(G)$ has rank $3$, then we must have
$[G,G]=G^{p}$.
Note that $|G^p|=p^2$, so suppose by contradiction, that $|G'|=p$. Observe that $G^p\le Z(G)$, so $G/Z(G)$ is a vector space over $\F_p$. Then, the commutator map in $G$ induces a non-degenerate alternating form on $G/Z(G)$, and
so $\dim_{\F_p}(G/Z(G))$ is even. This is a contradiction since $G/Z(G)$ has rank $3$.
We have thus shown that $[G,G]=G^{p}$. In order to distinguish further between different cases, we next turn our attention to $[\Omega_{1}(G),G]$. Notice that $\Omega_{1}(G)=\langle a\rangle G^{p}$. 
As $a\not\in Z(G)$ either $|[\Omega_{1}(G),G]|$  is of size $p$ or $p^{2}$. \\ \\
{\it Case 3}. ($Z(G)^{p}=1$, $G/Z(G)$ of rank $3$ and $|[\Omega_{1}(G),G]|=p$). Without loss of generality we can assume that $[\Omega_{1}(G),G]=\langle c^{p}\rangle$. There are two possibilities. Either $c\in C_{G}\left(\Omega_1\left(G\right)\right)=C_{G}(a)$ or not.
Suppose first that 
$c\in C_G\left(\Omega_{1}(G)\right)$. Then we have $[a,c]=1$, and we can pick $b$ such that $[a,b]=c^{p}$.
Replacing $b$ by $bc^{l}$ does not change these relations and thus we can assume that $[b,c]=b^{p\alpha}$ for some $0<\alpha<p$. If we let $\beta$ be the inverse of $\alpha$ modulo $p$ and we replace $a,c$ by $a^{\beta},c^{\beta}$, then we arrive at a group with presentation
         $$G_{18}=\langle a,b,c:\,a^{p}=b^{p^2}=c^{p^{2}}=1,\, [a,b]=c^{p},\,[b,c]=b^{p}\rangle.$$ 
Notice that this is a powerfully solvable group with a powerfully abelian chain $G>\langle b,c\rangle>\langle b\rangle>1$.
Suppose now $c\not\in C_G\left(\Omega_1(G)\right)$. Since $|[\Omega_1(G),G]|=|[a,G]|=p$, it follows that the conjugacy class of $a$ has order $p$, and so $|G:C_G(a)|=p$. Thus, we can pick $b$ such that $b\in C_G(a)$ and $[a,b]=1$. Replacing $a$ by a suitable power of $a$ we can suppose that $[a,c]=c^p$. As before, replacing $b$ by $bc^{l}$ does not change these relations, so we can also assume $[b,c]=b^{\alpha p}$ for some $0<\alpha<p$. Finally, if we let $\beta$ be the inverse of $\alpha$ modulo $p$ and we replace $c$ by $c^{\beta}$, we arrive at a group with presentation
$$
G_{19}=\langle a,b,c:\,a^{p}=b^{p^2}=c^{p^{2}}=1,\,[a,c]=c^p,\,[b,c]=b^{p}\rangle.
$$
This group is powerfully solvable with powerfully abelian chain $G>\langle b,c\rangle >\langle b\rangle >1$.
\\ \\
{\it Case 4}. ($Z(G)^{p}=1$, $G/Z(G)$ of rank $3$ and $|[\Omega_{1}(G),G]|=p^2$). In this case, commutation with $a$ induces a bijective linear map
\begin{eqnarray*}    
        F_{a}:G/\Omega_{1}(G) &\longrightarrow & G^{p} \\ 
           x\Omega_{1}(G) &\longmapsto & [a,x].
\end{eqnarray*}
Identifying $x\Omega_{1}(G)$ with $x^{p}$, we can think of $F_{a}$ as a linear operator on a two dimensional vector space over $\F_p$. Also replacing $b,c$ by a suitable $ba^{r},ca^{s}$ we can assume throughout that $[b,c]=1$. All the groups are going to be powerfully solvable with powerfully abelian chain $G>\langle b,c\rangle  >1$. 
\\ \\
{\it Case 4.1}. ($F_{a}$ is a scalar multiplication). Notice that this property still holds if we replace $a$ by any power of $a$ and thus it is independent of what $a$ we pick in $\Omega_{1}(G)\setminus G^{p}$. This is thus a characteristic property of $G$. Replacing $a$ with 
a power of itself we can assume that $F_{a}$ is the identity map. This gives us the group
  $$G_{20}=\langle a,b,c:\, a^{p}=b^{p^{2}}=c^{p^{2}}=1,\, [a,b]=b^{p},\,[a,c]=c^{p}\rangle.$$
{\it Case 4.2}. ($F_{a}$ is not a scalar multiplication). Again we see that this is a characteristic property of $G$. We can now pick $b$ and $c$ such that 
                     $$[a,b]=c^{p},\ [a,c]=b ^{p\alpha}c^{p\beta}.$$
Notice that the matrix for $F_{a}$ is 
                $$\left[\begin{array}{cc}
                           0  &  \alpha \\
 \mbox{}                          1  & \beta 
                  \end{array}\right]$$
with determinant $-\alpha$. This is an invariant for the given $a$ that does not depend on our choice of $b$ and  $c$. If we replace $a$ by $a^{r}$ and $c$ by $c^r$ then we get
              $$[a,b]=c^p,\ [a,c]=b^{p\alpha r^2}c^{p\beta r},$$
and the new determinant becomes $-\alpha r^{2}$. Pick some fixed $\tau$ such that $-\tau$ is a non-square in $\F_p$.
With appropriate choice of $r$ we can then assume that the determinant of $F_{a}$ is $-\alpha$ where either $\alpha=-1$
or $\alpha=\tau$. We thus have a group with one of the two presentations 
       $$G_{21}(\beta)=\langle a,b,c:\,a^{p}=b^{p^{2}}=c^{p^{2}}=1,\ [a,b]=c^{p},\ [a,c]=b^{-1}c^{p\beta},\,[b,c]=1\rangle,$$
and
     $$G_{22}(\beta)=\langle a,b,c:\,a^{p}=b^{p^{2}}=c^{p^{2}}=1,\ [a,b]=c^{p},\ [a,c]=b^{\tau}c^{p\beta},\,[b,c]=1\rangle.$$
Suppose we pick a different $\bar{b}=b^{r}c^{s}$. Then for $\alpha\in \{-1,\tau\}$ we have 
      $$[a,\bar{b}]=[a,b]^{r}[a,c]^{s}=c^{pr}(b^{p{\alpha}}c^{p\beta})^{s}=b^{ps\alpha}c^{p(r+s\beta)}=\bar{c}^{p}$$
where $\bar{c}=b^{s\alpha}c^{r+s\beta}$. Then 
\begin{eqnarray*}
      [a,\bar{c}] & = & [a,b]^{s\alpha}[a,c]^{r+s\beta} \\
                     & = &  c^{ps\alpha}(b^{p \alpha}c^{p\beta})^{r+s\beta} \\
                    & = & (b^{r}c^{s})^{p\alpha}\cdot (b^{s\alpha}c^{r+s\beta})^{p\beta} \\
                   & = & \bar{b}^{p\alpha}\bar{c}^{p\beta}.
\end{eqnarray*}
This shows that for the given $\alpha \in \{-1,\tau\}$, the constant $\beta\in \F_p$ is an invariant and we get $p$ distinct groups $G_{21}(\beta)$ and $p$ distinct groups $G_{22}(\beta)$. \\ \\
Adding up we have $7$ abelian groups and the groups $G_{6},\ldots , G_{20}, G_{21}(\beta), G_{22}(\beta)$, giving us in total
$\boldsymbol{22+2p}$ groups of order $p^{5}$. \\ \\
 Notice that we have seen that all powerful groups of order up to and including $p^{5}$ are powerfully solvable. Now take a powerful group of order $p^{6}$. Suppose it has a generator $a$ of order $p$, say $G=\langle a,H\rangle$
where $H<G$. Notice that $H$ is then powerful of order $p^{5}$ and thus powerfully solvable. As $[G,G]\leq H^{p}$ we then
see that $G$ is powerfully solvable. Thus all powerful groups of order $p^{6}$ are powerfully solvable with the possible exceptions of some groups of type $(2,2,2)$. We will see later that there are a number of groups of type $(2,2,2)$ that are not powerfully solvable. 
\mbox{}

\section{Growth}
Let $G$ be a powerfully solvable group of order $p^{n}$. From Theorem \ref{theorem presentation} and the discussion in Section 3, we know that we may assume that $G=\langle a_{1},\ldots ,a_{y},
a_{y+1},\ldots ,a_{y+x}\rangle$ where $o(a_{1})=\cdots =o(a_{y})=p$ and $o(a_{y+1})=\cdots =o(a_{y+x})=p^{2}$. Furthermore the generators can be chosen such that $|G|=p^{y+2x}$ and
         $$[a_{j},a_{i}]=a_{i+1}^{p\alpha_{i+1}(i,j)}\cdots a_{y+x}^{p\alpha_{y+x}(i,j)},$$
for $1\leq i<j\leq y+x$, where $0\leq \alpha_{k}(i,j)\leq p-1$ for $k=i+1,\ldots ,y+x$. For each such pair $(i,j)$ where $1\leq i\leq y$ there are $p^{x}$ possible relations for $[a_{j},a_{i}]$. There are $yx+{y\choose 2}$ such pairs. On the other hand, for a pair $(i,j)$
where $y+1\leq i\leq y+x$, for each given $i$ there are $y+x-i$ such pairs and $p^{y+x-i}$  possible relations $[a_{j},a_{i}]$. Adding up we see that the number of solvable presentations is $p^{h(x)}$ where 
 \begin{eqnarray*}
    h(x) & = & \left(yx+{y\choose 2}\right)x+1^{2}+2^{2}+\cdots (x-1)^{2} \\
          & = & \left(
          (n-2x)x+{n-2x\choose 2}\right)x+\frac{x(2x-1)(x-1)}{6} \\
          & = & \frac{1}{3}x^{3}-\frac{(2n-1)}{2}x^{2}+\frac{3n(n-1)+1}{6}x.
\end{eqnarray*} 
Thus 
      $$
      h'(x)=x^{2}-(2n-1)x+\frac{3n(n-1)+1}{6},
      $$
whose roots are $\frac{2n-1}{2}-\sqrt{\frac{1}{2}n^{2}-\frac{n}{2}+\frac{1}{12}}$ and $\frac{2n-1}{2}+\sqrt{\frac{1}{2}n^{2}-\frac{n}{2}+\frac{1}{12}}$. For large values of $n$ we have that the first root is between $0$ and $n/2$ whereas the latter is greater than $n$. 
Thus, for large $ n$,  the largest value of $h$ in the interval between $0$ and $n/2$ is $h(x(n))$ where $x(n)=\frac{2n-1}{2}-\sqrt{\frac{1}{2}n^{2}-\frac{n}{2}+\frac{1}{12}}$. Now $\lim_{n\rightarrow \infty}x(n)/n=1-\frac{1}{\sqrt{2}}$. Therefore
          $$\lim_{n\rightarrow \infty}\frac{h(x(n)}{n^{3}}=\lim_{n\rightarrow \infty}\frac{1}{3}(x(n)/n)^{3}-(x(n)/n)^{2}+\frac{1}{2}(x(n)/n))=\frac{-1+\sqrt{2}}{6}.$$
We now argue in a similar way as in [5]. Let $n$ be fixed. For any integer $x$ where $0\leq x\leq n/2$, let ${\mathcal P}(n,x)$ be the collection of all powerfully solvable presentations as above. It is not difficult to see that those presentations are consistent and thus the resulting group is of order $p^{n}$ and rank $n-x$. Furthermore $a_{1}^{p}=\cdots =a_{n-2x}^{p}=1$ and $a_{n-2x+1}^{p^{2}}=
\cdots =a_{n-x}^{p^{2}}=1$. We have just seen that, for large values of $n$, if we pick $x(n)$ such that the number of presentations is maximal then
                                 $$|{\mathcal P}(n,x(n))|=p^{\alpha n^{3}+o(n^{3})}$$
where $\alpha=\frac{-1+\sqrt{2}}{6}$. Let ${\mathcal P}_{n}$ be the total number of the powerfully solvable presentations where
$0\leq x\leq n/2$. Then ${\mathcal P}_{n}={\mathcal P}(n,0)\cup {\mathcal P}(n,1)\cup \cdots \cup {\mathcal P}(n,\lfloor n/2\rfloor)$ and thus
$$p^{\alpha n^{3}+o(n^{3})}=|{\mathcal P}(n,x(n))|\leq |{\mathcal P}_{n}|=|{\mathcal P}(n,0)|+\cdots +|{\mathcal P}(n,\lfloor n/2\rfloor)|\leq
n|{\mathcal P}(n,x(n))|=p^{\alpha n^{3}+o(n^{3})}.$$
This shows that $|{\mathcal P}_{n}|=p^{\alpha n^{3}+o(n^{3})}$. Let us show that this is also the growth of powerfully
solvable groups of exponent $p^{2}$ with respect to the order $p^{n}$. Clearly $p^{\alpha n^{3}+o(n^{3})}$ gives us an upper bound. We want to show that this is also a lower bound. Let $x=x(n)$ be as above and let $a_{1},\ldots ,a_{n-x}$ be a set of generators
for a powerfully solvable group $G$ where $a_{1}^{p}=\cdots =a_{n-2x}^{p}=1$ and $a_{n-2x+1}^{p^{2}}=
\cdots =a_{n-x}^{p^{2}}=1$. Notice that $\langle a_{1},\ldots ,a_{n-2x}\rangle G^{p}=\Omega_1(G)$, which is a  characteristic subgroup of $G$. It will be useful to consider a larger class of presentations for powerfully solvable groups of order $p^{n}$ where we still require $a_{1}^{p}=\cdots =a_{n-2x}^{p}=1$ and $a_{n-2x+1}^{p^{2}}=\cdots =a_{n-x}^{p^{2}}=1$.
We let ${\mathcal Q}(n,x)={\mathcal Q}(n,x(n))$ be the collection of all presentations with additional commutator relations 
                  $$[a_{i},a_{j}]=a_{1}^{p\alpha_{1}(i,j)}\cdots a_{n-x}^{p\alpha_{n-x}(i,j)}.$$
The presentation is included in ${\mathcal Q}(n,x)$ provided the resulting group is powerfully solvable of order $p^{n}$. Notice
that $G^{p}\leq Z(G)$ and as a result the commutator relations above only depend on the cosets $\overline{a_{1}}=a_{1}G^{p},\ldots , \overline{a_{n-x}}=a_{n-x}G^{p}$ and not on the exact values of $a_{1},\ldots ,a_{n-x}$. Consider the vector space
$V=G/G^{p}$ over $\F_p$ and let $W={\mathbb F}_{p}\overline{a_{1}}+\cdots +{\mathbb F}_{p}\overline{a_{n-2x}}$. Then let
                     $$H=\{\phi \in \mbox{GL}(n-x,p):\,\phi(W)=W\}.$$
There is now a natural action from $H$ on  ${\mathcal Q}(n,x)$. Suppose we have some presentation with generators $a_{1},\ldots ,a_{n-x}$ as above. Let $\phi\in H$ and suppose 
                      $$\overline{a_{i}}^{\phi}=\beta_{1}(i)\overline{a_{1}}+\cdots +\beta_{n-x}(i)\overline{a_{n-x}}.$$
We then get a new presentation in ${\mathcal Q}(n,x)$ for $G$ with respect to the generators $b_{1},\ldots ,b_{n-x}$ where
$b_{i}=a_{1}^{\beta_{1}(i)}\cdots a_{n-x}^{\beta_{n-x}(i)}$.  \\ \\
Suppose there are $l$ powerfully solvable groups of exponent $p^{2}$ and order $p^{n}$ where furthermore $|G^{p}|=p^{x}$. Pick powerfully solvable presentations $p_{1},\ldots ,p_{l}\in {\mathcal P}(n,x)$ for these. Let $q$ be powerfully solvable presentation in ${\mathcal P}(n,x)$ of a group $K$ with generators $b_{1},\ldots ,b_{n-x}$. Then $q$ is also a presentation for an isomorphic group $G$
with presentation $p_{i}$ and generators $a_{1},\ldots ,a_{n-x}$. Let $\phi:K\rightarrow G$ be an isomorphism and let $\psi:K/K^{p}\rightarrow G/G^{p}$ be the corresponding linear isomorphism. This gives us a linear automorphism $\tau\in H$ induced by $\tau(\overline{a_{i}})=\psi(\overline{b_{i}})$. Thus $q=p_{i}^{\tau}$. Therefore 
                $${\mathcal P}(n,x)\subseteq p_{1}^{H}\cup p_{2}^{H}\cup\cdots \cup p_{l}^{H}.$$
From this we get 
     $$p^{\alpha n^{3}+o(n^{3})}=|{\mathcal P}(n,x)|\leq |p_{1}^{H}|+\cdots +|p_{l}^{H}|\leq lp^{n^{2}},$$
and it follows that $l\geq p^{\alpha n^{3}+o(n^{3})}$. We thus get the following result.
\begin{theo} The number of powerfully solvable groups of exponent $p^{2}$ and order $p^{n}$ is 
$p^{\alpha n^{3}+o(n^{3})}$, where $\alpha=\frac{-1+\sqrt{2}}{6}$.
\end{theo}
\noindent As mentioned in [5] the growth of all powerful $p$-groups of exponent $p^{2}$ and order $p^{n}$ is $p^{\frac{2}{27}n^{3}+o(n^{3})}$. This claim was though not proved and we will fill in the details here. \\ \\
As before we consider a group $G$ of order $p^{n}=p^{y+2x}$ with generators $a_{1},\ldots ,a_{y+x}$ where $o(a_{1})=\cdots =o(a_{y})=p$ and $o(a_{y+1})=\cdots =o(a_{y+x})=p^{2}$. This time we can though include all powerful relations 
         $$[a_{j},a_{i}]=a_{1}^{p\alpha_{y+1}(i,j)}\cdots a_{y+x}^{p\alpha_{y+x}(i,j)}$$
for $1\leq i<j\leq y+x$, where $0\leq \alpha_{k}(i,j)\leq p-1$ for $k=y+1,\ldots ,y+x$. For each such pair $(i,j)$ there are $p^{x}$ possible relations for $[a_{j},a_{i}]$. We thus see that the number of  presentations is $p^{h(x)}$ where 
    $$h(x) = {y+x\choose 2}x =  {n-x\choose 2}x= \frac{x^{3}}{2}-\frac{(2n-1)}{2}x^{2}+\frac{n(n-1)}{2}x.$$
Thus 
      $$h'(x)=\frac{3}{2}\left(x^{2}-\frac{2(2n-1)}{3}x+\frac{n(n-1)}{3}\right)$$
and using the same kind of analysis as before we see that for a large $n$, $h$ takes its maximal value for $x(n)=\frac{2n-1}{3}-
\sqrt{\frac{n^{2}}{9}-\frac{n}{9}+\frac{1}{9}}$. Notice that $\lim_{n\rightarrow \infty}\frac{x(n)}{n}=1/3$. Therefore 
       $$\lim_{n\rightarrow \infty}\frac{h(x(n)}{n^{3}}
       =\lim_{n\rightarrow \infty}\frac{1}{2}
       \cdot\left(\frac{n-x(n)}{n}\right)
       \cdot\left(\frac{n-1-x(n)}{n}\right)
       \cdot\frac{x(n)}{n}
       =2/27.$$
The same argument as above shows then that the growth of all powerful groups of exponent $p^{2}$ with respect to order $p^{n}$ is 
$p^{\frac{2}{27}n^{3}+o(n^{3})}$. \\ \\
Later on we will be working with a special subclass ${\mathcal P}$ of powerful $p$-groups, namely those that are of type {$(2,\stackrel{r}{\ldots} ,2)$ with $r\ge 1$}. In this case the number of presentations for groups of order $p^{n}$, $n$ even, is $p^{h(n)}$ where $h(n)={n/2 \choose 2}n/2$ and
         $$\lim_{n\rightarrow \infty}\frac{h(n)}{n^{3}}=\lim_{n\rightarrow \infty}\frac{n/2(n/2-1)n/2}{2n^{3}}=1/16.$$
Thus the growth here is $p^{\frac{1}{16}n^{3}+o(n^{3})}$. 
\section{Groups of type $(2,\ldots,2)$}
We have seen that powerful nilpotence and powerful solvability is preserved under taking quotients. These properties however work badly under taking subgroups. Our next two results underscore this.
\begin{prop} Let $G$ be any powerful $p$-group of exponent $p^{2}$. There exists a powerfully nilpotent group $H$ of exponent $p^{2}$ and powerful class $2$ such that $G$ is powerfully embedded in $H$.
\end{prop}
\noindent{\bf Proof}\ \ Suppose $G=\langle a_{1},\ldots ,a_{r}\rangle$ where $a_{1}^{p}=\cdots =a_{s}^{p}=1$, $a_{s+1}^{p^{2}}=
\cdots =a_{r}^{p^{2}}=1$ and where $|G|=p^{s+2(r-s)}$. Let $N=\langle x_{s+1}\rangle \times \cdots\times \langle x_{r}\rangle$
be a direct product of cyclic groups of order $p^{2}$. Let  $H=(G\times N)/M$, where $M=\langle a_{s+1}^{p}x_{s+1}^{-p},
\ldots ,a_{r}^{p}x_{r}^{-p}\rangle$. Notice that $[G,H]=[G,G]\leq G^{p}$ and thus $G$ is powerfully embedded in $H$. Also, as $[H,H]=G^{p}=N^{p}$,
we see that 
             $$1\leq \langle x_{1},\ldots ,x_{r}\rangle\leq H$$
is a powerfully central chain and thus $H$ is powerfully nilpotent of powerful class at most $2$. $\Box$ \\ \\
{\bf Remark}. (1) There exist powerful $p$-groups of exponent $p^{2}$ that are not powerfully solvable and thus a powerfully embedded subgroup of a powerfully nilpotent group of powerful class $2$ does not
even need to be powerfully solvable. \\
(2) There exist powerfully nilpotent groups of exponent $p^{2}$ that are of arbitrary large powerful class and so the proposition above shows that a powerfully nilpotent group of powerful class $2$ could have 
a powerfully embedded powerfully nilpotent subgroup of arbitrary large powerful class. \\ \\
Next result shows that the subgroup structure of a powerfully nilpotent group of powerful class $2$ is even more arbitrary. Notice that such a group is in particular nilpotent of class $2$ and it turns out that any finite $p$-group
of class $2$ can occur as a subgroup. 
\begin{prop} Let $G$ be any finite  $p$-group of nilpotency class $2$. There exists a powerfully nilpotent group $H$ of powerful class $2$ that contains $G$ as a subgroup.
\end{prop}
\noindent{\bf Proof}\ \ Suppose $[G,G]$ has a basis $a_{1},\ldots ,a_{m}$ as an abelian group where $o(a_{i})=p^{j_i}$. Let $N=\langle x_{1}\rangle\times\cdots \times \langle x_{m}\rangle$ be a direct product
of cyclic groups where $o(x_{i})=p^{j_i+1}$. Now let $H=(G\times N)/M$ where $M=\langle a_{1}x_{1}^{-p},\ldots, a_{m}x_{m}^{-p}\rangle$. Then $G$ embeds as a subgroup
of $H$. Notice also that 
       $$1\leq \langle x_{1},\ldots ,x_{m}\rangle \leq H$$
is powerfully central and thus $H$ is powerfully nilpotent of powerful class $2$. $\Box$ \\ \\ 
Thus powerful nilpotence and powerful solvability are in general not as satisfactory as notions for powerful groups as nilpotence
and solvability for the class of all groups. For a rich subclass of powerful groups things however turn out much better. This is the class ${\mathcal P}$ of all powerful groups of type $(2,\stackrel{r}{\ldots} ,2)$ that we considered in Section 5. \\ \\
For a group $G\in {\mathcal P}$ we have that $G^{p}\leq Z(G)$. It follows that the map $G/G^{p}\rightarrow G^{p},\,aG^{p}\mapsto a^{p}$ is a bijection and therefore, for any $H\geq G^{p}$, we have $|H/G^{p}|=|H^{p}|$. 
\begin{lemm} Let $G\in {\mathcal P}$ and $H,K\leq G$ where $G^{p}\leq K$. Then $H^{p}\cap K^{p}=(H\cap K)^{p}$.
\end{lemm}
\noindent{\bf Proof}\ \ We have 
\begin{eqnarray*}
    |H^{p}\cap K^{p}|
      =   |(HG^{p})^{p}\cap K^{p}|
    & = & \frac{|(HG^{p})^{p}|\cdot|K^{p}|}{|(HK)^{p}|}\\
    & = & \frac{|HG^{p}/G^{p}|\cdot |K/G^{p}|}{|HK/G^{p}|}\\
    & = & |(HG^{p}\cap K)/G^{p}| \\
    & = & |(H\cap K)G^{p}/G^{p}|\\
    & = & |(H\cap K)^{p}|.
\end{eqnarray*}        
As $(H\cap K)^{p}\leq H^{p}\cap K^{p}$ it follows that $H^{p}\cap K^{p}=(H\cap K)^{p}$. $\Box$ 
\begin{theo} Let $G$ be a powerfully nilpotent group in ${\mathcal P}$ and let $H$ be a powerful subgroup of $G$. Then $H$ is
powerfully nilpotent of powerful class less than or equal to the powerful class of $G$.
\end{theo}
\noindent{\bf Proof}\ \ Suppose $G$ has powerful nilpotence class $c$ and that we have a powerfully central chain $G=G_{0}>G_{1}>\cdots >G_{c}=1$. As $G^{p}\leq Z(G)$ and $(G^{p})^{p}=1$, multiplying a term by $G^{p}$ makes no difference. Also as the 
powerful class is $c$ we get a strictly decreasing powerfully central chain $G=G_{0}>G_{1}G^{p}>\cdots >G_{c-1}G^{p}>1$. Without loss of generality we can thus assume that $G_{1},\ldots ,G_{c-1}$ contain $G^{p}$ as a subgroup. We claim that 
        $$H=H\cap G_{0}\geq H\cap G_{1}\geq \cdots \geq H\cap G_{c-1}\geq 1$$
is powerfully central. Using Lemma 6.3 we have
$$
[H\cap G_{i},H]\leq [H,H]\cap [G_{i},G]\leq H^{p}\cap G_{i+1}^{p}=(H\cap G_{i+1})^{p},
$$
for $0\leq i\leq c-1$. Hence $H$ is powerfully nilpotent of powerful class at most $c$. $\Box$ \\ \\ 
\begin{theo} Let $G$ be a powerfully solvable group in ${\mathcal P}$ and let $H$ be a powerful subgroup of $G$. Then $H$ is powerfully solvable of powerful derived length less than or equal to the powerful derived length of
$G$.
\end{theo}
\noindent{\bf Proof}\ \ Suppose the powerful derived length of $G$ is $d$ and that we have a powerfully abelian chain $G=G_{0}>G_{1}>\cdots >G_{d}=1$. Arguing like in the proof of the previous theorem, we can assume that $G_{d-1}$ contains $G^{p}$. We show that 
                   $$H=H\cap G_{0}\geq H\cap G_{1}\geq \cdots \geq H\cap G_{c-1}\geq H\cap G_{c}=1$$
is a powerfully abelian chain. Using Lemma 6.3, we have  $[H\cap G_{i},H\cap G_{i}]\leq [H,H]\cap [G_{i},G_{i}]\leq H^{p}\cap G_{i+1}^{p}=(H\cap G_{i+1})^{p}$. This shows that $H$ is powerfully solvable of powerful derived length at most $d$. $\Box$ \\ \\ 
We introduce some useful notation. We use $H\leq_{\mathcal P} G$ to stand for $H,G\in {\mathcal P}$ and $H\leq G$.  We use $H\unlhd_{\mathcal P}G$ for $H,G\in {\mathcal P}$ and
$H$ powerfully embedded in $G$. The notations  $H<_{\mathcal P}G$ and $H\lhd _{\mathcal P}G$ are defined naturally in a similar way. \\ \\
Let $G$ be a powerful $p$-group in $\mathcal{P}$ and let $V=G/G^{p}$ be the associated vector space over $\F_p$.
The structure of $G$ is determined by the commutator relations 
 \begin{equation}
                  [a,b]=c^{p},
\end{equation}
where there exists such $c\in G$ for each pair $a,b$ in $G$.  Notice that $[a,b]$ and $c^{p}$ only depend on the cosets $aG^{p},bG^{p}$ and $cG^{p}$. Identifying
the two vector spaces $G/G^{p}$ and $G^{p}$ under the map $G/G^{p}\rightarrow G^{p},xG^{p}\mapsto x^{p}$, we get a natural alternating product on $V$ with 
the relations (1) translating to 
                        $$[aG^{p},bG^{p}]=cG^{p}.$$                     
Let ${\mathcal G}$ be the collection of all powerful subgroups of $G$ that are of type $(2,\ldots ,2)$ and let ${\mathcal V}$ be the collection of all the alternating
subalgebras of $V$. So for $U$ to be a subalgebra of $V$ it needs to be a subspace where $[U,U]\leq U$. Notice that $[H,H]\leq H^{p}$ translates to
$[HG^{p}/G^{p},HG^{p}/G^{p}]\leq HG^{p}/G^{p}$. Recall that for $H,K\in {\mathcal G}$ we write $H\trianglelefteq_{\mathcal P} K$ for $H$ powerfully embedded in $K$. For $U,W\in {\mathcal V}$ we likewise
write $U\trianglelefteq W$ for $U$ an ideal of $W$. Notice that $[H,K]\leq H^{p}$ translates to $[HG^{p}/G^{p},KG^{p}/G^{p}]\leq HG^{p}/G^{p}$. \\ \\
If $H,G\in{\mathcal P}$ such that $G$ is powerfully nilpotent and $H\trianglelefteq_{\mathcal{P}} G$, then the quotient $G/H$ has naturally the structure
of  powerful group of type $(2,\ldots ,2)$ with $[aH,bH]=[a,b]H$. \\ \\
{\bf Definition}. We say that a group $G\in\mathcal{P}$ is {\it powerfully simple} if $G\not =1$ and if $H\lhd_{\mathcal P} G$ implies that $H=1$. \\ \\
{\bf Definition}. Let $H,G\in\mathcal{P}$ with $H\lhd_{\mathcal P}G$. We say that $H$ is 
a maximal powerfully embedded ${\mathcal P}$-subgroup of $G$ if there is no $H<K<G$ such that $K\unlhd_{\mathcal P}G$. 
\begin{lemm}
Let $G\in\mathcal{P}$. We have that $H$ is a maximal powerfully embedded $\mathcal{P}$-subgroup of $G$ if and only if $G/H$ is powerfully simple.
\end{lemm}
\noindent{\bf Proof}\ \ Let $H<K<G$. Now as $H$ is powerful of type $(2,\ldots,2)$ we have $H\cap G^{p}=H^{p}$. Therefore  
$[K,G]\leq K^{p}H$ if and only if
$$[K,G]\leq (K^{p}H)\cap G^{p}=K^{p}(H\cap G^{p})=K^{p}H^{p}=K^{p}.$$
The result follows from this. $\Box$ \\ \\
{\bf Remark}. Let $H,K\in {\mathcal G}$ and let $U$ and $W$ be the associated alternating algebras in ${\mathcal V}$. Suppose that $H$ is powerfully
embedded in $K$. Then $K/H$ is powerfully simple if and only if $W/U$ is a simple alternating algebra and the latter happens if and only if $U$ is
a maximal ideal of $V$. \\ \\
We will next prove a Jordan-Holder type theorem for alternating algebras. Suppose $A,B,a,b\in {\mathcal V}$ where $A\lhd B$ and $a\lhd b$. Let 
${\mathcal I}_{A}^{B}=\{Z:A\leq Z\leq B\}$ and ${\mathcal I}_{a}^{b}=\{x:a\leq x\leq b\}$. We get natural projections $P:{\mathcal I}_{a}^{b}\rightarrow
{\mathcal I}_{A}^{B}$ and $Q:{\mathcal I}_{A}^{B}\rightarrow {\mathcal I}_{a}^{b}$ given by
      $$P(x)=A+B\cap x\mbox{\  and\ }Q(Z)=a+b\cap Z.$$
\begin{lemm}We have $P(a)\trianglelefteq P(b)$ and $Q(A)\trianglelefteq Q(B)$.  Furthermore $P(b)/P(a)$ is isomorphic to $Q(B)/Q(A)$.
\end{lemm}
\noindent{\bf Proof}\ \ Notice that $P(a)=A+B\cap a, P(b)=A+B\cap b, Q(A)=a+b\cap A$ and $Q(B)=a+b\cap B$.  As $A\trianglelefteq B$, we have $[P(a),P(b)]=[A+B\cap a,A+B\cap b]\leq A+[B\cap a,B\cap b]$. Now as $B$ is a subalgebra and $a\trianglelefteq b$ we have that
this is contained in $A+B\cap a$. The second claim follows from this by symmetry. \\ \\
Now for $P(b)/P(a)$, notice first that we have
$$
B\cap b\cap (A+B\cap a)=B\cap b\cap A+B\cap a=A\cap b+B\cap a,
$$
and for $u,v,w\in B\cap b$ it follows that
  $$       [u,v]+A+B\cap a=w+A+B\cap a \Leftrightarrow [u,v]+A\cap b+B\cap a=w+A\cap b+B\cap a.$$
By symmetry 
          $$[u,v]+a+b\cap A=w+a+b\cap A \Leftrightarrow [u,v]+a\cap B+b\cap A=w+A\cap b+B\cap a.$$
The isomorphism of $P(b)/P(a)$ and $P(B)/B(A)$ follows from this. $\Box$ \\ \\
The Jordan-Holder theorem for alternating algebras is proved from this in the standard way. \\ \\
{\bf Definition}. Let $V$ be an alternating algebra. A chain $0=U_{0}\lhd U_{1} \lhd \cdots \lhd U_{n}=V$ is a {\it composition series} for $V$ if all the factors $U_{1}/U_{0},\ldots ,U_{n}/U_{n-1}$ are simple alternating algebras.
\begin{theo} Let $V$ be an alternating algebra. Then all composition series have the same length and same composition factors up to order.
\end{theo}
\noindent{\bf Definition}. Let $G\in {\mathcal P}$. A chain $1=H_{0}\lhd_{\mathcal P} H_{1}\lhd_{\mathcal P}\cdots \lhd_{\mathcal P} H_{n}=G$ is a {\it powerful composition series} for $G$ if all the factors $H_{1}/H_{0},
\ldots ,H_{n}/H_{n-1}$ are powerfully simple.
\begin{theo}Let $G$ be a powerful $p$-group of type $(2,\ldots ,2)$ with two powerful composition series, say
     $$1=H_{0}\lhd_{\mathcal P}H_{1}\lhd_{\mathcal P}\cdots \lhd_{\mathcal P}H_{n}=G$$
and
    $$1=K_{0}\lhd_{\mathcal P}K_{1}\lhd_{\mathcal P}\cdots \lhd_{\mathcal P}K_{m}=G.$$
Then $m=n$ and the powerfully simple factors $H_{1}/H_{0},H_{2}/H_{1},\ldots ,H_{n}/H_{n-1}$ are isomorphic to $K_{1}/K_{0}$, $K_{2}/K_{1}$, $\ldots ,K_{n}/K_{n-1}$
(in some order). 
\end{theo}
\noindent{\bf Proof}\ \ Replace the terms $H_{i},K_{j}$ by their associated alternating algebras $U_{i}, V_{j}$. The result now follows from the Jordan-Holder theorem for alternating algebras. $\Box$ \\ \\
{\bf Definition}. We refer to the unique factors of a powerful composition series of a group $G\in {\mathcal P}$ as the {\it powerful composition factors} of $G$. 
\begin{coro}
A group $G\in {\mathcal P}$ is powerfully solvable if and only if the powerful composition factors are cyclic of order $p^{2}$.
\end{coro}
\noindent{\bf Proof}\ \ Any powerful abelian chain of $G$ can be refined to a chain with factors that are cyclic of order $p^{2}$. $\Box$

\section{The classification of powerfully simple groups of type $(2,2,2)$.}
From previous section we know that this task is equivalent to classifying all simple alternating algebras of dimension $3$. \\ \\
Following [4], any given alternating algebra of dimension $3$
over ${\mathbb F}_p$ can be represented by a $3\times 3$ matrix over ${\mathbb F}_p$. Here the matrix 
$$\left[\begin{array}{lll}
   \alpha_{11} & \alpha_{12} & \alpha_{13} \\
  \alpha_{21} & \alpha_{22} & \alpha_{23} \\
  \alpha_{31} & \alpha_{32} & \alpha_{33}
\end{array}\right]$$
corresponds to the $3$-dimensional alternating algebra ${\mathbb F}_{p}v_{1}+{\mathbb F}_{p}v_{2}+{\mathbb F}_{p}v_{3}$ where
\begin{eqnarray*}
    v_{2}v_{3} & = & \alpha_{11}v_{1}+\alpha_{21}v_{2}+\alpha_{31}v_{3} \\
     v_{3}v_{1} & = & \alpha_{12}v_{1}+\alpha_{22}v_{2}+\alpha_{32}v_{3} \\
   v_{1}v_{2} & = & \alpha_{13}v_{1}+\alpha_{23}v_{2}+\alpha_{33}v_{3}.
\end{eqnarray*}
From last section we know that this corresponds to a powerful $p$-group of order $p^{6}$ with generators $a_{1},a_{2},a_{3}$ of order $p^{2}$ satisfying the relations:
\begin{eqnarray*}
    [a_{2},a_{3}] & = & a_{1}^{p\alpha_{11}}a_{2}^{p\alpha_{21}}a_{3}^{p\alpha_{31}} \\
 \mbox{}    [a_{3},a_{1}] & = & a_{1}^{p\alpha_{12}}a_{2}^{p\alpha_{22}}a_{3}^{p\alpha_{32}} \\
 \mbox{}  [a_{1},a_{2}] & = & a_{1}^{p\alpha_{13}}a_{2}^{p\alpha_{23}}a_{3}^{p\alpha_{33}}.
\end{eqnarray*}
In [4] it is shown that two such matrices $A,B$ represent the
same alternating algebra with respect to different basis if and only there exists an invertible $3\times 3$ matrix $P$ such that
        $$B=\frac{1}{\mbox{det\,}(P)}P^{t}AP.$$
We write $B\simeq  A$ if they are related in this way. This turns out to be slightly more general than being congruent (that is $B=P^{t}AP$).
\begin{lemm} Let $\lambda \in {\mathbb F}_p^{*}$. Then $\lambda A\simeq A$.
\end{lemm} 
\noindent{\bf Proof}\ \ Let $P=\frac{1}{\lambda}I$. Then $\frac{1}{\det(P)}P^{t}AP=\lambda^{3}\cdot \frac{1}{\lambda^{2}}A=\lambda A$. $\Box$ \\ \\
From this we easily get the following corollary.
\begin{prop}
We have $B\simeq A$ if and only if there exists $C$ such that $A$ is congruent to $C$ and $B=\lambda C$. 
\end{prop}
\noindent In particular two matrices that are congruent are equivalent. We can write each such matrix $A$ in a unique way as a sum of a symmetric and an anti-symmetric matrix, namely $A_{s}=\frac{A+A^{t}}{2}$ and $A_{a}=\frac{A-A^{t}}{2}$. As shown in [4] we have that 
$A\simeq B$ if and only if $A_{s}\simeq B_{s}$ and $A_{a}\simeq B_{a}$. We will determine all the equivalence classes and therefore all 
powerful $p$-groups of type $(2,2,2)$. From this we will then single out those that are powerfully simple. \\ \\
{\bf Classification of the symmetric matrices}.
It is known that every symmetric matrix is congruent to a diagonal matrix and furthermore to exactly one  of the following $D(1,1,1)$, $ D(\tau,1,1)$, $D(1,1,0)$, $ D(\tau, 1, 0)$, $D(1,0,0)$, $D(\tau,0,0)$ and $D(0,0,0)$, where $\tau$ is a fixed non-square in ${\mathbb F}_p^{*}$ and $D(\alpha,\beta,\gamma)$ is the $3\times 3$ matrix with $\alpha, \beta$ and $\gamma$ on the diagonal (compare [1, Chapter 6, Theorem 2.7]). Now $D(1,1,1)$ is equivalent to $\tau D(1,1,1)$ where the latter has determinant $\tau$ modulo $({\mathbb F}_p^{*})^{2}$. Hence $D(1,1,1)$ and $D(\tau,1,1)$ are equivalent. Also $D(\tau,0,0)=\tau D(1,0,0)$ is equivalent $D(1,0,0)$. When the rank is $2$ then multiplying the matrix by a constant $\lambda\in {\mathbb F}_p^{*}$ doesn't change the value of the determinant modulo $({\mathbb F}_p^{*})^{2}$. Hence $D(1,1,0)$ and
$D(\tau,1,0)$ are not equivalent. Up to equivalence we thus get only $5$  matrices:
      $$D(1,1,1),\, D(1,1,0),\, D(\tau,1,0),\, D(1,0,0)\, \mbox{and }D(0,0,0).$$ \\
{\bf Classification of the anti-symmetric matrices.}
The situation regarding the anti-symmetric matrices is simpler as there are only two equivalence classes. One containing the zero matrix and one for the non-zero matrices. This comes from the fact that there are only two alternating forms (up to isomorphism) for a $3$-dimensional algebra $V$. Either $V^{\perp}$ has dimension  $3$ or $1$. \\ \\
{\bf Classification of the alternating algebras.}
Let $A$ be some $3\times 3$ matrix over ${\mathbb F}_p$ and let $V$ be the corresponding alternating algebra. The symmetric part of $A$ equips $V$ with a corresponding symmetric bilinear form $\langle \ ,\ \rangle_{s}$ and the anti-symmetric part of $A$ equips 
$V$ with a corresponding alternating form $\langle \ ,\ \rangle_{a}$. Now there are two possibilities for $\langle \ ,\ \rangle_{a}$. If it is zero then $A$ is symmetric and we get ${\bf five}$ alternating algebras corresponding to the $5$ diagonal matrices listed above. From now
on we can thus assume that $\langle \ ,\ \rangle_{a}$ is non-zero. Thus $V^{\perp_{a}}$ is of dimension $1$. Say 
        $$V^{\perp_{a}}={\mathbb F}_pv_{3}$$
so that
$$
V=({\mathbb F}_{p}v_{1}+{\mathbb F}_{p}v_{2})\operp_{a} {\mathbb F}_{p}v_{3}
$$
for some $v_1,v_2\in V$. For our classification we will divide first into $3$ cases. For Case 1, we have $\langle v_{3},v_{3}\rangle_{s}\not =0$. For Case 2, we have $\langle v_{3},v_{3}\rangle_{s}=0$ and $(V^{\perp_{a}})^{\perp_{s}}=({\mathbb F}_pv_{3})^{\perp_{s}}=V$. Finally
for Case 3, we have $\langle v_{3},v_{3}\rangle_{s}=0$ and $(V^{\perp_{a}})^{\perp_{s}}=({\mathbb F}_pv_{3})^{\perp_{s}}<V$. \\ \\
{\it Case 1}. 
We can here find a basis $v_{1},v_{2},v_{3}$ for $V$ where 
            $$V={\mathbb F}_{p}v_{1}\operp_{s} {\mathbb F}_{p}v_{2}\operp_{s} {\mathbb F}_{p}v_{3}.$$
%\\ \\
%
{\it Case 1.1}.  Suppose first that the rank of $\langle\ ,\ \rangle_{s}$ is $1$. In this case it is easy to see that we can pick our basis further so that
$$\begin{array}{lll}
    \langle v_{1},v_{2}\rangle_{a}=1,  & \langle v_{1},v_{3}\rangle_{a}=0, & \langle v_{2},v_{3}\rangle_{a}=0, \\
   \langle v_{1},v_{1}\rangle_{s}=0,  & \langle v_{2},v_{2}\rangle_{s}=0, & \langle v_{3},v_{3}\rangle_{s}=1.
\end{array}$$ 
Indeed, notice that we can always multiply the relevant matrix by a constant to get $\langle v_{3},v_{3}\rangle_{s}=1$ and then it is easy to pick our $v_{1},v_{2}$ such that $\langle v_{1},v_{2}\rangle_{a}=1$. In this case we thus have only $\boldsymbol{1}$ algebra. \\ \\
{\it Case 1.2}. Suppose next that the rank of $\langle\ ,\ \rangle_{s}$ is $2$.  Here again by multiplying by a constant we can assume that $\langle v_{3},v_{3}\rangle_{s}=1$ and we can assume that $\langle v_{2},v_{2}\rangle_{s}=0$. Now $\langle v_{1},v_{1}\rangle_{s}=\lambda^{2}$ or $\langle v_{1},v_{1}\rangle_{s}=\tau\lambda^{2}$ for some $\lambda\in {\mathbb F}_p^{*}$. By replacing $v_{1}$ by $\frac{1}{\lambda}v_{1}$ we can assume that $\langle v_{1},v_{1}\rangle_{s}$ is either $1$ or $\tau$. Notice that we have also seen above that these cases are genuinely distinct. Now that $v_{1}$ has been chosen we can replace $v_{2}$ by a suitable multiple to ensure that $\langle v_{1},v_{2}\rangle_{a}=1$. We thus get $\boldsymbol{2}$ algebras
 $$\begin{array}{lll}
    \langle v_{1},v_{2}\rangle_{a}=1,  & \langle v_{1},v_{3}\rangle_{a}=0, & \langle v_{2},v_{3}\rangle_{a}=0, \\
   \langle v_{1},v_{1}\rangle_{s}=1,  & \langle v_{2},v_{2}\rangle_{s}=0, & \langle v_{3},v_{3}\rangle_{s}=1.
\end{array}$$ 
and
$$\begin{array}{lll}
    \langle v_{1},v_{2}\rangle_{a}=1,  & \langle v_{1},v_{3}\rangle_{a}=0, & \langle v_{2},v_{3}\rangle_{a}=0, \\
   \langle v_{1},v_{1}\rangle_{s}=\tau,  & \langle v_{2},v_{2}\rangle_{s}=0, & \langle v_{3},v_{3}\rangle_{s}=1.
\end{array}$$  \\ \\
{\it Case 1.3}. We are then only left with the case where the rank of $\langle \ ,\ \rangle_{s}$ is $3$. It is not difficult to see that in this case we can pick $v_{1},v_{2},v_{3}$ such that 
$$\begin{array}{lll}
    \langle v_{1},v_{2}\rangle_{a}=1,  & \langle v_{1},v_{3}\rangle_{a}=0, & \langle v_{2},v_{3}\rangle_{a}=0, \\
   \langle v_{1},v_{1}\rangle_{s}=\alpha,  & \langle v_{2},v_{2}\rangle_{s}=1, & \langle v_{3},v_{3}\rangle_{s}=1,
\end{array}$$ 
where $\alpha\in {\mathbb F}_p^{*}$. We want to see when we get an equivalent algebra by changing $\alpha$ to $\beta$. If we multiply the presentation by a constant it must be by a square if we still want $\langle v_{3},v_{3}\rangle_{s}=1$. Say, we multiply 
by $\lambda^{2}$ and then replace $v_{3}$ by $\frac{1}{\lambda} v_{3}$. Notice that we now have 
                $$\langle v_{1},v_{2}\rangle_{a}=\lambda^{2},\ \langle v_{1},v_{1}\rangle_{s} =\alpha\lambda^{2},\ \langle v_{2},v_{2}\rangle_{s}=\lambda^{2}.$$
We are now looking for all possible $\bar{v}_{1}=av_{1}+bv_{2}$ and $\bar{v}_{2}=cv_{1}+dv_{2}$ where $\langle \bar{v}_{1},\bar{v}_{2}\rangle_{a}=1$, $\langle \bar{v}_{1},\bar{v}_{2}\rangle_{s}=0$ and  $\langle \bar{v}_{2},\bar{v}_{2}\rangle_{s}=1$. This gives us the following system of equations:
\begin{eqnarray*}
      \lambda^{2}(ad-bc) & = & 1 \\
    \lambda^{2}(\alpha ac+bd) & = & 0 \\
\lambda^{2}(\alpha c^{2}+d^{2}) & = & 1.
\end{eqnarray*}
We look first for all the solutions where $c=0$. Notice that in this case we must have  $\lambda^{2} ad=1$, $\lambda^2d^2=1$ and $b=0$. Thus $\langle \bar{v}_{1},\bar{v}_{1}\rangle_{s}=\lambda^{2}(\alpha a^{2}+b^{2})=\lambda^{2}\frac{\alpha}{d^{2}\lambda^{4}}=\frac{\alpha}{\lambda^{2}d^{2}}=\alpha$. \\ \\
Next we look for solutions where $c\not =0$ but $d=0$. Then we must have $\lambda^{2}bc=-1$, $\lambda^2\alpha c^2=1$ and $a=0$. Here $\langle \bar{v}_{1},\bar{v}_{1}\rangle_{s}=\lambda^{2}(\alpha a^{2}+b^{2})= \frac{\lambda^{2}}{c^{2}\lambda^{4}}=\frac{\alpha}{\alpha c^{2}\lambda^{2}}=\alpha$. \\ \\
Finally we are left with finding all solutions where $cd\not =0$. Then $a=-\frac{bd}{\alpha c}$, and the first equation above gives us
               $$1=-\lambda^{2}\left(\frac{bd^{2}}{\alpha  c}+bc\right)=-\frac{b}{\alpha c}\cdot \lambda^{2}(d^{2}+\alpha c^{2})=-\frac{b}{\alpha c}.$$
Thus $b=-\alpha c$ and $a=-\frac{bd}{\alpha c}=d$. Therefore 
\begin{eqnarray*}
       \langle \bar{v}_{1},\bar{v}_{1}\rangle & = & \lambda^{2}(\alpha a^{2}+b^{2}) \\
                              & = & \lambda^{2}(\alpha d^{2}+\alpha^{2}c^{2}) \\
                              & = & \alpha \lambda^{2}(\alpha c^{2}+d^{2}) \\
                             & = & \alpha.
\end{eqnarray*}
We have thus seen that the value of $\alpha$ doesn't change and we have $\boldsymbol{p-1}$ different algebras here.  \\ \\
{\it Case 2}.
Here we are assuming that $\langle v_{3},v_{3}\rangle_{s}=0$ and that $v_{3}$ is orthogonal to $v_{1},v_{2}$ as well. Again we consider few subcases. \\ \\
{\it Case 2.1}. Suppose that the rank of $\langle \ ,\ \rangle_{s}$ is zero. Then clearly we have $\boldsymbol{1}$ algebra. 
$$\begin{array}{lll}
    \langle v_{1},v_{2}\rangle_{a}=1,  & \langle v_{1},v_{3}\rangle_{a}=0, & \langle v_{2},v_{3}\rangle_{a}=0, \\
   \langle v_{1},v_{1}\rangle_{s}=0,  & \langle v_{2},v_{2}\rangle_{s}=0, & \langle v_{3},v_{3}\rangle_{s}=0.
\end{array}$$  \\ \\
{\it Case 2.2}. Suppose next that  the rank of $\langle \ ,\ \rangle_{s}$ is $1$. By multiplying by a suitable constant we can assume that $\langle v_{1},v_{1}\rangle_{s}=1$ and $\langle v_{2},v_{2}\rangle_{s}=0$. Finally replacing $v_{2}$ by an appropriate 
multiple we can  also assume that $\langle v_{1},v_{2}\rangle_{a}=1$. We thus also get here only $\boldsymbol{1}$ algebra
$$\begin{array}{lll}
    \langle v_{1},v_{2}\rangle_{a}=1,  & \langle v_{1},v_{3}\rangle_{a}=0, & \langle v_{2},v_{3}\rangle_{a}=0, \\
   \langle v_{1},v_{1}\rangle_{s}=1,  & \langle v_{2},v_{2}\rangle_{s}=0, & \langle v_{3},v_{3}\rangle_{s}=0.
\end{array}$$  \\ \\
{\it Case 2.3}. Finally we are left with the case when the rank of $\langle \ ,\ \rangle_{s}$ is $2$. Here it is easy to see that we can pick our basis further so that 
$$\begin{array}{lll}
    \langle v_{1},v_{2}\rangle_{a}=1,  & \langle v_{1},v_{3}\rangle_{a}=0, & \langle v_{2},v_{3}\rangle_{a}=0, \\
   \langle v_{1},v_{1}\rangle_{s}=\alpha ,  & \langle v_{2},v_{2}\rangle_{s}=1, & \langle v_{3},v_{3}\rangle_{s}=0.
\end{array}$$ 
Similar calculations as for Case 1.3 show that we get distinct algebras for different values  of $\alpha$. Thus here we have $\boldsymbol{p-1}$ algebras. \\ \\
{\it Case 3}.
Here we are assuming that $\langle v_{3},v_{3}\rangle_{s}=0$ but that $v_{3}$ is not orthogonal to everything in $V$ with respect to $\langle\ ,\ \rangle_s$. Thus $({\mathbb F}_pv_{3})^{\perp_{s}}$ has dimension $2$. Suppose
                   $$({\mathbb F}_pv_{3})^{\perp_{s}}={\mathbb F}_pv_{2}+{\mathbb F}_pv_{3}.$$
It is not difficult to see that we can pick our basis such that
$$\begin{array}{llll}
    \langle v_{1},v_{2}\rangle_{a}=1,  & \langle v_{1},v_{3}\rangle_{a}=0, & \langle v_{2},v_{3}\rangle_{a}=0, &\\
   \langle v_{1},v_{1}\rangle_{s}=0, & \langle v_{1},v_{2}\rangle_{s}=0,  & \langle v_{1},v_{3}\rangle_{s}=1, & \langle v_{2},v_{3}\rangle_{s}=0.
\end{array}$$ 
Now there are two subcases.  \\ \\
{\it Case 3.1}. If the rank of $\langle\ ,\ \rangle_{s}$ is $2$ then we must have $\langle v_{2},v_{2}\rangle_{s}=0$ and this gives us $\boldsymbol{1}$ algebra. \\ \\
{\it Case 3.2}. If the rank of $\langle \ , \ \rangle_{s}$ is $3$ then $\langle v_{2},v_{2}\rangle_{s}\not =0$ and after multiplying by a suitable constant we can assume that this value is $1$ (and then afterwards adjust things so that the other assumptions hold again).
Thus we get again $\boldsymbol{1}$ algebra.  \\ \\
Adding up we see that in total we get $12+2(p-1)$ algebras and thus the same number of powerful $p$-groups of type $(2,2,2)$.
Before listing these we state and prove a proposition that shows how we can see which of these are powerfully simple.
\begin{prop} An alternating algebra $V$ over ${\mathbb F}_{p}$ of dimension $3$ is simple if and only if $V\cdot V=V$. 
\end{prop}
\noindent{\bf Proof}\ \ This condition is clearly necessary as $V\cdot V$ is an ideal of $V$. To see that it is sufficient, suppose $V\cdot V=V$ and let
$I$ be a proper ideal. We want to show that $I=0$. We argue by contradiction and suppose $I$ is an ideal of dimension either $1$ or $2$. If $I$ is of dimension $2$, then $V/I$ is $1$ dimensional and thus we get the contradiction that $V\cdot V\leq I<V$. Now suppose $I$ is of dimension $1$, say $V=I+{\mathbb F}_{p}v_{1}+{\mathbb F}_{p}v_{2}$. Then $V\cdot V\leq I+{\mathbb F}_{p}v_{1}v_{2}$. But the dimension of $I+{\mathbb F}v_{1}v_{2}$ is at most $2$ and we get the contradiction that $V\cdot V<V$. $\Box$ \\ \\
We have thus determined the presentation matrices up to equivalence and got in total $12+2(p-1)$. As we described at the beginning of the section this gives us a classification of all the alternating algebras of dimension $3$ over ${\mathbb F}_{p}$ that in turn gives us a classification of all the powerful $p$-groups of type $(2,2,2)$. Furthermore, last proposition tells us how we read from the presentation whether a given alternating algebra is simple and thus whether the corresponding powerful group is powerfully simple. 
The work above gives us the following list of powerful $p$-groups of type $(2,2,2)$. As the power relations for all of these are
$a_{1}^{p^{2}}=a_{2}^{p^{2}}=a_{3}^{p^{2}}=1$ we omit these below. Here $\tau$ is a fixed non-square in ${\mathbb F}_p$. \\
$$\begin{array}{l}
   A_{1}=\langle a_{1},a_{2},a_{3}:\,[a_{2},a_{3}]=a_{1}^{p},\,[a_{3},a_{1}]=a_{2}^{p},\,[a_{1},a_{2}]=a_{3}^{p}\rangle; \\
 A_{2}=\langle a_{1},a_{2},a_{3}:\,[a_{2},a_{3}]=a_{2}^{-p},\,[a_{3},a_{1}]=a_{1}^{p},\,[a_{1},a_{2}]=a_{3}^{p}\rangle; \\
 A_{3}=\langle a_{1},a_{2},a_{3}:\,[a_{2},a_{3}]=a_{1}^{p}a_{2}^{-p},\,[a_{3},a_{1}]=a_{1}^{p},\,[a_{1},a_{2}]=a_{3}^{p}\rangle; \\
 A_{4}=\langle a_{1},a_{2},a_{3}:\,[a_{2},a_{3}]=a_{1}^{p\tau}a_{2}^{-p},\,[a_{3},a_{1}]=a_{1}^{p},\,[a_{1},a_{2}]=a_{3}^{p}\rangle; \\
 A_{5}=\langle a_{1},a_{2},a_{3}:\,[a_{2},a_{3}]=a_{2}^{-p}a_{3}^{p},\,[a_{3},a_{1}]=a_{1}^{p}a_{2}^{p},\,[a_{1},a_{2}]=a_{1}^{p}\rangle; \\
 A_{6}(\alpha)=\langle a_{1},a_{2},a_{3}:\,[a_{2},a_{3}]=a_{1}^{p\alpha}a_{2}^{-p},\,[a_{3},a_{1}]=a_{1}^{p}a_{2}^{p},\,[a_{1},a_{2}]=a_{3}^{p}\rangle,\ 1\leq\alpha\leq p-2; \\
 B_{1}=\langle a_{1},a_{2},a_{3}:\,[a_{2},a_{3}]=a_{1}^{-p}a_{2}^{-p},\,[a_{3},a_{1}]=a_{1}^{p}a_{2}^{p},\,[a_{1},a_{2}]=a_{3}^{p}\rangle; \\
 B_{2}=\langle a_{1},a_{2},a_{3}:\,[a_{2},a_{3}]=a_{1}^{p},\,[a_{3},a_{1}]=a_{2}^{p},\,[a_{1},a_{2}]=1\rangle; \\
 B_{3}=\langle a_{1},a_{2},a_{3}:\,[a_{2},a_{3}]=a_{1}^{p\tau},\,[a_{3},a_{1}]=a_{2}^{p},\,[a_{1},a_{2}]=1\rangle; \\
 B_{4}=\langle a_{1},a_{2},a_{3}:\,[a_{2},a_{3}]=a_{2}^{-p},\,[a_{3},a_{1}]=a_{1}^{p},\,[a_{1},a_{2}]=1\rangle; \\
 B_{5}=\langle a_{1},a_{2},a_{3}:\,[a_{2},a_{3}]=a_{1}^{p}a_{2}^{-p},\,[a_{3},a_{1}]=a_{1}^{p},\,[a_{1},a_{2}]=1\rangle; \\
 B_{6}=\langle a_{1},a_{2},a_{3}:\,[a_{2},a_{3}]=a_{2}^{-p}a_{3}^{p},\,[a_{3},a_{1}]=a_{1}^{p},\,[a_{1},a_{2}]=a_{1}^{p}\rangle \\
 B_{7}(\alpha)=\langle a_{1},a_{2},a_{3}:\,[a_{2},a_{3}]=a_{1}^{p\alpha}a_{2}^{-p},\,[a_{3},a_{1}]=a_{1}^{p}a_{2}^{p},\,[a_{1},a_{2}]=1\rangle,\ 1\leq \alpha\leq p-2;\\
 C_{1}=\langle a_{1},a_{2},a_{3}:\,[a_{2},a_{3}]=a_{1}^{-p}a_{2}^{-p},\,[a_{3},a_{1}]=a_{1}^{p}a_{2}^{p},\,[a_{1},a_{2}]=1\rangle; \\
 C_{2}=\langle a_{1},a_{2},a_{3}:\,[a_{2},a_{3}]=a_{1}^{p},\,[a_{3},a_{1}]=1,\,[a_{1},a_{2}]=1\rangle; \\
 D=\langle a_{1},a_{2},a_{3}:\,[a_{2},a_{3}]=1,\,[a_{3},a_{1}]=1,\,[a_{1},a_{2}]=1\rangle.
\end{array}$$
\mbox{}\\
Of these $12+2(p-1)$ groups, the groups $A_{1},\ldots ,A_{6}(\alpha)$ are the powerfully simple groups. There are $5+(p-2)$ of these.

\end{document}